\documentclass[12pt]{amsart}
\usepackage{amscd}
\def\frk{\frak}               

\def\mm{{\frk m}}

\def\Phi{{\frk n}}
\def\Phi{{\frk N}}
\def\opn#1#2{\def#1{\operatorname{#2}}} 
\opn\chara{char}
\opn\length{\ell}
\opn\pd{pd}
\opn\rk{rk}
\opn\projdim{proj\,dim}
\opn\injdim{inj\,dim}
\opn\rank{rank}
\opn\depth{depth}
\opn\grade{grade}
\opn\height{height}
\opn\embdim{emb\,dim}
\opn\codim{codim}

\opn\Tr{Tr}
\opn\bigrank{big\,rank}
\opn\superheight{superheight}\opn\lcm{lcm}
\opn\trdeg{tr\,deg}
\opn\reg{reg}
\opn\lreg{lreg}
\opn\ini{in}
\opn\lpd{lpd}
\opn\size{size}
\opn\div{div}
\opn\Div{Div}
\opn\cl{cl}
\opn\Cl{Cl}
\opn\Spec{Spec}
\opn\Supp{Supp}
\opn\supp{supp}
\opn\Sing{Sing}
\opn\Ass{Ass}
\opn\Min{Min}
\opn\Ann{Ann}
\opn\Rad{Rad}
\opn\Soc{Soc}
\opn\Im{Im}
\opn\Ker{Ker}
\opn\Coker{Coker}
\opn\Am{Am}
\opn\Hom{Hom}
\opn\Tor{Tor}
\opn\Ext{Ext}
\opn\End{End}
\opn\Aut{Aut}
\opn\id{id}

\opn\nat{nat}
\opn\pff{pf}
\opn\Pf{Pf}
\opn\GL{GL}
\opn\SL{SL}
\opn\mod{mod}
\opn\ord{ord}
\opn\Gin{Gin}
\opn\Hilb{Hilb}
\opn\aff{aff}
\opn\con{conv}
\opn\relint{relint}
\opn\st{st}
\opn\lk{lk}
\opn\cn{cn}
\opn\core{core}
\opn\vol{vol}
\opn\link{link}
\opn\star{star}
\opn\gr{gr}

\def\pot#1#2{#1[\kern-0.28ex[#2]\kern-0.28ex]}

\opn\dirlim{\underrightarrow{\lim}}
\opn\inivlim{\underleftarrow{\lim}}
\let\iso=\cong

\let\to=\rightarrow

\def\Implies{\ifmmode\Longrightarrow \else
        \unskip${}\Longrightarrow{}$\ignorespaces\fi}
\def\implies{\ifmmode\Rightarrow \else
        \unskip${}\Rightarrow{}$\ignorespaces\fi}
\def\iff{\ifmmode\Longleftrightarrow \else
        \unskip${}\Longleftrightarrow{}$\ignorespaces\fi}

\let\:=\colon
\newtheorem{Theorem}{Theorem}[section]
\newtheorem{Lemma}[Theorem]{Lemma}
\newtheorem{Corollary}[Theorem]{Corollary}
\newtheorem{Proposition}[Theorem]{Proposition}
\newtheorem{Remark}[Theorem]{Remark}

\newtheorem{Example}[Theorem]{Example}

\let\epsilon\varepsilon
\let\phi=\varphi
\let\kappa=\varkappa
\def\qed{\ifhmode\textqed\fi
      \ifmmode\ifinner\quad\qedsymbol\else\dispqed\fi\fi}
\def\textqed{\unskip\nobreak\penalty50
       \hskip2em\hbox{}\nobreak\hfil\qedsymbol
       \parfillskip=0pt \finalhyphendemerits=0}
\def\dispqed{\rlap{\qquad\qedsymbol}}

\opn\dis{dis}
\def\pnt{{\raise0.5mm\hbox{\large\bf.}}}

\opn\Lex{Lex}

\begin{document}

\title{Monomial cycle basis on Koszul homology modules}

\author{  Dorin Popescu}
\subjclass{13D02, 13D07, 13P10, 13D99}
\thanks{The author was mainly supported by Marie Curie Intra-European
Fellowships MEIF-CT-2003-501046 and partially supported by   the
Ceres program 4-131/2004 of the Romanian Ministery of Education and
Research}

\address{Dorin Popescu, Institute of Mathematics "Simion Stoilow", University of Bucharest,
P.O.Box 1-764, Bucharest 014700, Romania}
\email{dorin.popescu@imar.ro} \maketitle

\begin{abstract}
It gives a class of $p$-Borel principal ideals of a polynomial
algebra over a field $K$ for which the graded Betti numbers do not
depend on the characteristic of $K$ and the Koszul homology modules
have monomial cyclic basis. Also it shows that all principal
$p$-Borel ideals have binomial cycle basis on Koszul homology
modules.
\end{abstract}

\section*{Introduction}

Let $K$ be an infinite  field, $S=K[x_1,\ldots,x_n]$, $n\geq 2$ the
polynomial ring over a field $K$ and $I\subset S$ a graded ideal.
Consider the reverse lexicographical order on the monomials of $S$.
 Let $M$ be a graded $S$-module and
$\beta_{ij}(M)=\beta_{ij}$  the graded Betti numbers of $M$. The
Castelnuovo-Mumford regularity of $M$ is
$\reg(M)=\max\{j-i:\beta_{ij}(M)\not =0\}$. By a theorem of Bayer
and Stillman  \cite{BS} we have $\reg(\Gin(I))=\reg(I)$.  If char
$K=0$ then $\Gin(I)$ is strongly stable, that is, it is monomial and
for each monomial $u$ of $\Gin(I)$ and $1\leq j<i\leq n$ such that
$x_i|u$ it follows $x_j(u/x_i)\in \Gin(I)$. Then $\reg(\Gin(I)$ is
the highest degree of  minimal generators of $\Gin(I)$ by Eliahou
and Kervaire \cite{EK}. If char $K=p>0$ then Borel fixed ideals are
just the so called $p$-Borel ideals and they are not necessarily
strongly stable and it is hard to give a formula for the regularity
of these ideals. Let $I$ be a monomial ideal of $S$, $u$ a monomial
of $I$ and $\nu_i(u)$ be the highest power of $x_i$ dividing $u$.
Let $a,b$ be two integers and $a=\Sigma_{i\geq 0}a_i p^i$,
$b=\Sigma_{i\geq 0}b_i p^i$ be the $p$-adic expansion of $a$,
respectively $b$. We say that $a\leq_p b$ if $a_i\leq b_i$ for all
$i$.
 It is well known that a monomial
ideal $I$ is  $p$-Borel  if for any monomial $u\in I$ and $1\leq
j<i\leq n$ and a positive integer $t$ such that $t\leq_p\nu_i(u)$ it
holds $x_j^t(u/x_i^t)\in I$. This is a pure combinatorial
description of the $p$-Borel ideals which can be given independently
of the characteristic of $K$. Let $u$ be a monomial of $S$ and
$J=\langle u\rangle$  the smallest monomial ideal containing $u$.
$J$ is called principal $p$-Borel ideal. For such ideals there
exists a complicated formula for regularity in terms of $u$
conjectured by Pardue \cite{Pa} and proved in two papers \cite{AH},
\cite{HP} (another proof is given in \cite{HPV}).

 In general it
is hard to bound the regularity of a graded ideal $I$. If char $K=0$
and $d(I)$ is the highest degree of a minimal monomial generator of
$I$ then D. Bayer and D. Mumford \cite{BM} showed that $\reg(I)\leq
(2d(I))^{2^{n-2}}$. Caviglia and Sbarra \cite{CS} showed that the
same bound holds for all characteristic of $K$. This bound seems to
be sharp since  Mayr and Meyer \cite{MM} gave an example with
$d(I)=4$ and $\reg(I)\geq 2^{2^{n-1}}+1$. Thus in general a bound
for the regularity  is very high, but what about if we restrict to
some classes of ideals? In \cite{Po} it is showed that  $d(I)\leq
\reg(I)\leq nd(I)$ if $I$ is a $p$-Borel ideal.  The proof uses the
formula conjectured by Pardue. The Betti numbers and the regularity
of a graded ideal $I$ can depend of characteristic of the field $K$
even when $I$ is monomial. This is not the case when $I$ is strongly
stable as follows from Eliahou and Kervaire resolution \cite{EK}.
Using again the formula conjectured by Pardue we get that the
regularity of principal $p$-Borel ideals does not depend on the
characteristic of $K$ (see \cite{Po}).

Let $Syz_t(M)=\Ker(F_t\to F_{t-1})$ be  the $t$-th syzygy module of
$M$. The module $M$ is called {\em $(r,t)$-regular}  if $Syz_t(M)$
is $(r+t)$-regular in the sense that all generators of $F_j$ for
$t\leq j\leq s$ have degrees $\leq j+r$. The {\em $t$-regularity} of
$M$, that is, the regularity of $Syz_t(M)$ is given by
$$(t-\reg)(M)=\min\{r\ : M\;  \text{is} \;  (r,t)-\text{regular}\}.$$
Obviously, we have $(t-\reg)(M)\leq ((t-1)-\reg)(M)$. If the
inequality is strict and $r=(t-\reg)(M)$, then $(t,r)$ is called a
{\em corner} of $M$ and $\beta_{t,r+t}(M)$ is an {\em extremal Betti
number} of $M$ \cite{BCP}. As the regularity of principal $p$-Borel
ideals is completely determined by their extremal Betti numbers it
is natural to ask if these numbers depend on the characteristic of
$K$. They do not depend indeed in a more general frame explained
bellow.

Bayer and Stillman  proved that if $I$ is Borel-fixed then it
satisfies the following property:
$$(I:x_j^{\infty})=(I:(x_1,\ldots,x_j)^{\infty })$$
\noindent for $j=1,\ldots,n$. A monomial ideal $I\subset S$
satisfying the above condition is said   to be of {\sl Borel type}
\cite{HPV}. For a monomial $u$, let $\nu_i(u)$ be the highest power
of $x_i$ which divides  $u$ and $m(u)=\max\{i: \nu_i(u)\not = 0\}$.
For a monomial ideal $I\not =0$, let $G(I)$ be the unique set of
monomial minimal generators of $I$ and $m(I)=\max\{m(u): u\in
G(I)\}$. We define recursively an ascending chain of monomial
ideals:
$$I=I_0\subset I_1 \subset \ldots $$
\noindent as follows: We let $I_0=I$ and if $I_e\not =S$ is already
defined then set $I_{e+1}=(I_e:x_{n_e}^{\infty})$ for $n_e=m(I_e)$,
$n_0>n_1>\ldots >n_e>\ldots$. If $I_e=S$ then the chain ends. Let
$(I_e)_{0\leq e\leq q}$, $(n_e)$ be the sequences obtained above,
$I_q\not =0$ and $I_{q+1}=S$. Let $S_i=K[x_1,\ldots,x_{n_i}]$,
$J_i=I_i\cap S_i$ and $J_i^{sat}$ the saturation of $J_i$ made in
$S_i$.

\begin{Theorem}[\cite{Po}]
 \label{extrBorel} Let $I\subset S$ be a Borel type
ideal. Then $S/I$ has at most $q+1$-corners among
$(n_i,s(J_i^{sat}/J_i))$, $0\leq i\leq q$ and the corresponding
extremal Betti numbers are
$$\beta_{n_i,s(J_i^{sat}/J_i)+n_i}\ \
(S/I)=\dim_K(J_i^{sat}/J_i)_{s(J_i^{sat}/J_i)},$$ where
$s(N)=\max\{i:N_i\not =0\}$ for a graded $S$-module $N$ of finite
length. In particular the corners of $S/I$ and their corresponding
extremal Betti numbers do not depend on the characteristic of $K$.
Moreover the regularity of $S/I$ does not depend of characteristic
of $K$.
\end{Theorem}

Remains to ask what about the general Betti numbers? Bellow we
remind you a nice case when this is true.
 The principal $p$-Borel ideals $I\subset S$ such that $S/I$ is
Cohen-Macaulay have the form $I=\Pi_{j=0}^s(m^{[p^j]})^{\alpha_j}$,
$0\leq \alpha_j<p$, where $m^{[p^j]}=(x_1^{p^j},\ldots,x_n^{p^j})$.
For these ideals is well known the description of a canonical
monomial cycle basis of the Koszul homology module $H_i(x;S/I)$
given by Aramova and Herzog \cite{AH} (see details in \ref{ah}). One
can easily see from this description that  $\beta_{ij}(S/I)$ does
not depend on the characteristic of the field $K$ for all $i,j$.

Now  let $I$ be the $p$-Borel ideal generated by the monomial
$x_{n-1}^{\gamma}x_n^{\alpha}$ for some integer $\gamma,\alpha\geq
0$, that is
$$I=\Pi_{j=0}^s((m_{n-1}^{[p^j]})^{\gamma_j}(m^{[p^j]})^{\alpha_j}),$$
where $m_{n-1}=(x_1,\ldots,x_{n-1})$, and $\gamma_j,\alpha_j$ are
defined by the $p$-adic expansion of $\gamma$, respectively
$\alpha$.
 Suppose that $\alpha_j+\gamma_j<p$ for all $0\leq j\leq s$. Then
$H_i(x;S/I)$ has a monomial cycle basis for all $i\geq 2$, and
$\beta_{ij}(S/I)$ does not depend on the characteristic of $K$ for
all $i,j$ (see \ref{main}).

We saw that in some cases of principal $p$-Borel ideals there exist
a monomial cycle basis for the homology modules of $S/I$. How it is
in general? If $I\subset S$ is a monomial ideal then $H_2(x;S/I)$
has a monomial cycle basis (see \ref{2cyc}). Unfortunately, in
general there are no monomial cycle basis even on the Koszul
homology modules of principal $p$-Borel ideals as shows our Example
\ref{bi}. However if $I$ is a principal $p$-Borel ideal then
$H_3(x;S/I)$ has a binomial cycle basis (see \ref{main1}). In
general our Example \ref{five} shows that there are reduced monomial
ideals which have not even a trinomial cycle basis. Perhaps in
general there exist monomial reduced ideals $I$ in $n$-variables
such that there exists non-zero cycles of length $n-2$ which are not
modulo bounds sum of cycles of length $<n-2$.

We express our  thanks to J. Herzog especially for some
discussions around Theorem \ref{main} and Lemma \ref{n-1}.

\section{Cycles  of Koszul homology modules of
monomial ideals}

Let $S=K[x_1,\ldots,x_n]$ be a polynomial algebra over a field $K$
and $I\subset S$ a monomial ideal. A cycle $z\in K_i(x;S/I)$ has the
form $z=\Sigma_{j=1}^s\gamma_ju_j e_{\sigma_j}$, $\gamma_j\in K^*$,
$u_j$ monomials, $\sigma_j\subset \{1,\ldots n\}$, $|\sigma_j|=i$
for all $1\leq j\leq s$. Since $I$ is monomial the Koszul
antiderivation $\partial $ is multigraded and each cycle is a sum of
multigraded cycles. The cycle $z$ is multigraded if
$u_jx_{\sigma_j}=u_1x_{\sigma_1}$ for all $1\leq j\leq s$, here
$x_{\sigma_1}=\Pi_{k\in\sigma_1}x_k$. We denote
$m(u_j)=\max\{i;x_i|u_j\}$ and $m(\sigma_j)=m(x_{\sigma_j})$.
 Note that in $z$  we may
suppose $\sigma_j\not =\sigma_t$ for $j\not =t$ because otherwise it
follows $u_j=u_t$ ($z$ is multigraded) and so we may reduce the sum.
The element $u_je_{\sigma_j}$ is a monomial cycle if $\partial
(u_je_{\sigma_j})=0$, that is $x_tu_j\in I$ for all $t\in \sigma_j$.

We introduce a totally order on the monomial elements $ue_{\sigma}$
of $K_i(x;S/I)$ ($u$ monomial) by "$ue_{\sigma}\geq ve_{\tau}$" if
either "$u>_{rlex} v$" or "$u= v$" and "$x_{\sigma}\geq_{rlex}
x_{\tau}$", here rlex denotes the reverse lexicographical order on
the monomials of $S$. As usually we denote $in( z)=u_1e_{\sigma_1}$
if $u_1e_{\sigma_1}>u_je_{\sigma_j}$ for all $j>1$. A $\sigma_j$ is
called a {\sl neighbour} in $z$ of $\sigma_1$ if $|\sigma_j\setminus
\sigma_1|=1$.

\begin{Lemma}
\label{neighbour} If $\sigma_1$ has no neighbour in $z$ then
$u_1e_{\sigma_1}$ is a monomial cycle.
\end{Lemma}
\begin{proof}
Since $z$ is a cycle all the terms of $\partial (u_1e_{\sigma_1})$
should be reduced with terms of some $\partial (u_je_{\sigma_j})$,
$j>1$. But this is possible only if $\sigma_j$ is a neighbour of
$\sigma_1$.

\end{proof}

\begin{Lemma}
\label{bound} Let $z=\Sigma_{j=1}^s\gamma_ju_j e_{\sigma_j}$ be a
multigraded cycle, $\gamma_j\in K^*$, $u_j$ monomials,
$\sigma_j\subset \{1,\ldots n\}$, $|\sigma_j|=i$ for all $1\leq
j\leq s$. Then the following statements hold:

\begin{enumerate}
\item{} If $in (z)=u_1e_{\sigma_1}$ and $m(u_1)> m(\sigma_1)$ then
there exists a multigraded element $w\in B_i(x;S/I)$ such that $in
(w)=in( z)$.

\item{} For every multigraded cycle $w$ there exist a multigraded
cycle $z$ of the above form  in the same multigraded homology class
with $w$ such that $m(u_j)\leq m(\sigma_j)$ for all $1\leq j\leq s$.

\item{} If $z$ is in the form given by (2) it follows
$m(\sigma_j)=m(\sigma_1)$ for all $1\leq j\leq s$.
\end{enumerate}
\end{Lemma}
\begin{proof} (1) Take a $q>m(\sigma_1) $ such that $x_q|u_1$
and set $y=(u_1/x_q)e_{\sigma_1\cup \{q\}}$. We have $\partial y$ is
the sum of $(u_1/x_q)(\partial e_{\sigma_1}) \wedge e_{ \{q\}}$ with
+ or - $u_1e_{\sigma_1}$. Thus $in( \partial y)=in( z)$.

(2)+(3) Substracting from $z$ such elements $w$ of $B_i(x;S/I)$ we
may arrive to the case $m(u_j)\leq m(\sigma_j)$. Since $z$ is
multigraded we get then $m(\sigma_j)=m(\sigma_1)$.

\end{proof}

\begin{Lemma}
\label{xr}Let $z=\Sigma_{j=1}^s\gamma_ju_j e_{\sigma_j}$ be a
multigraded cycle as in the above Lemma. Suppose that $m(u_j)\leq
m(\sigma_j)$ for all $1\leq j\leq s$. Then $x_ru_j\in I$ for all
$1\leq j\leq s$.
\end{Lemma}
\begin{proof}
By Lemma \ref{bound} (3) we get $r=m(\sigma_1)=m(\sigma_j)$. The
terms $x_r u_j e_{\sigma_j\setminus\{r\}}$ of $\partial
(u_je_{\sigma_j})$ cannot be reduced since $\sigma_j\setminus\{r\}$
are all different. It follows necessarily $x_ru_j\in I$ since $z$ is
a cycle.
\end{proof}

 Let ${\mathcal M}_i(x;S/I)$ be the subspace of $K_i(x;S/I)$
generated by all monomial cycles.

\begin{Lemma}
\label{2-cycle} Let $z=\Sigma_{j=1}^s\gamma_ju_j e_{\sigma_j}$ be a
multigraded 2-cycle, $\gamma_j\in K^*$, $u_j$ monomials,
$\sigma_j\subset \{1,\ldots ,n\}$, $|\sigma_j|=i$ for all $1\leq
j\leq s$. Suppose that $m(u_j)\leq m(\sigma_j)$ for all $j$, $s>1$
and $in(z)=u_1e_{\sigma_1}$. Then one of the following conditions
holds:

\begin{enumerate}
\item{} $in( z)$ is a monomial cycle,
\item{} $in (z)\equiv u_je_{\sigma_j}$ mod $(B_2(x;S/I)+{\mathcal
M}_2(x;S/I))$ for some $1< j\leq s$.
\end{enumerate}
\end{Lemma}

\begin{proof} By Lemma \ref{bound} (3) and our hypothesis we get
$r=m(\sigma_1)=m(\sigma_j)$ for all $1\leq j\leq s$.  Let
$\sigma_1=\{a,r\}$. If $x_au_1\in I$ then $u_1e_{\sigma_1}$ is a
monomial cycle. Otherwise by Lemma \ref{neighbour} there exists a
neighbour $\sigma_j=(\sigma_1\setminus\{a\})\cup\{b\}$ for some
$1\leq b< r$. As $in(z)=u_1e_{\sigma_1}$ we have $\sigma_1>\sigma_j$
and so $a<b$. From $x_au_1=x_bu_j$ ($z$ is multigraded!) it follows
$x_b|u_1$. Set $y=(u_1/x_b) e_{\{a,b,r\}}\in K_3(x;S/I)$. We have
$$\partial
y=-u_1e_{\sigma_1}+u_je_{\sigma_j}+(x_ru_1/x_b)e_{\{a,b\}}.$$ Using
Lemma \ref{xr} we have $x_ru_t\in I$ for all $1\leq t\leq s$. This
shows that the last term of $\partial y$ is a monomial cycle, which
is enough.
\end{proof}

\begin{Theorem}
\label{2cyc} Every 2-cycles of $K_2(x;S/I)$ belongs to
$B_2(x;S/I)+{\mathcal M}_2(x;S/I)$, that is coincides modulo
$B_2(x;S/I)$ with a sum of monomial cycles. In particular,
$H_2(x;S/I)$ has a monomial cycle basis.
\end{Theorem}

 \begin{proof} Note that given a 2-cycle $z$ in the form from Lemma
\ref{bound} (2) $in( z)$ can be substitute in $z$ modulo
$B_2(x;S/I)$ with one monomial term smaller than $in( z)$ and some
monomial cycles which can be removed from $z$ (see \ref{2-cycle}).
By recurrence we arrive finally to the case when $z$ has just one
term which must be then monomial cycle.
\end{proof}

\section{Some useful examples}

The purpose of this paper is to study when the Koszul homology
modules of principal $p$-Borel ideals $I$ have monomial cycle basis.
This is the case when $I$ is the smallest $p$-Borel ideal containing
a power $u$ of one variable $x_r$ (the so  called $p$-Borel ideal
generated by $u$), that is $I=\Pi_{j\geq 0}
(\mm_r^{\alpha_j})^{[p^j]}$, where $\mm_r=(x_1,\ldots,x_r)$ and
$\alpha_j<p$ are non-negative integers (see \cite{AH}). For an ideal
$J$ we denote by $J^{[p^j]}$ the ideal generated in $S$ by
$\phi(J)$, $\phi$ being the $K$-automorphism of $S$ given by $x\to
x^{p^j}$. An interesting and promising example is the following:

\begin{Example}
\label{inter} {\em Let $n=4$,  $S=K[x_1,\ldots,x_4]$ and $I$ the
$p$-Borel ideal generated by the monomial $\{x_3x_4^p\}$, that is
$I=(x_1,x_2,x_3)(x_1^p,\ldots,x^p_4)$. Then
$z=x_1^{p-1}x_3x_4^{p-1}e_{124}-x_1^{p-1}x_2x_4^{p-1}e_{134}$
 is a cycle. Take $y=x_1^{p-1}x_4^{p-1}e_{1234}\in K_4(x;S/I)$.
 We have $\partial y
 =z+x_1^px_4^{p-1}e_{234}-x_1^{p-1}x_4^pe_{123}$. Note that
 $x_1^{p-1}x_4^p\in I$ but $x_1^px_4^{p-1}\not\in I$. Thus
 $z$ coincides with $x_1^px_4^{p-1}e_{234}$ modulo $B_3(x;S/I)$.
 As $z$ and $\partial y$ are cycles we get  $x_1^px_4^{p-1}e_{234}$
monomial cycle.}
\end{Example}

Unfortunately, in general there are not monomial cycle basis even on
the Koszul homology modules of principal $p$-Borel ideals as shows
the following examples:

\begin{Example}
\label{bi} {\em Let $n=4$, that is $S=K[x_1,\ldots,x_4]$ and $I$ the
$p$-Borel ideal generated by the monomial $\{x_2x_4^p\}$, that is
$I=(x_1,x_2)(x_1^p,\ldots,x_4^p)$. Consider the element
$z=x_2x_3^{p-1}x_4^{p-1}e_{134}-x_1x_3^{p-1}x_4^{p-1}e_{234}\in
K_3(x;S/I)$. We see that $z$ is binomial cycle but
$in(z)=x_2x_3^{p-1}x_4^{p-1}e_{134}$ is not a monomial cycle because
$x_1x_2x_3^{p-1}x_4^{p-1}\not\in I$. Note that $z$ is multigraded
and in its multigraded homology class take another element of the
form $z+\partial y$, $y\in K_4(x;S/I)$. Since $y$ must be
multigraded from the same multigraded class with $z$ we see that the
only possibility is to take $y=x_3^{p-1}x_4^{p-1}e_{1234}$. It
follows that there exist no monomial cycle in the homology class of
$z$. }
\end{Example}

\begin{Remark}
\label{tri}{\em Let  $S$ and $I$ as in Example \ref{bi}.
 The element
$$z=x_1^{p-1}x_3x_4^{p-1}e_{124}-x_1^{p-1}x_2x_4^{p-1}e_{134}+
x_1^px_4^{p-1}e_{234}\in K_3(x;S/I)$$
  is a cycle in the form given by Lemma \ref{bound} (2) but belong
  to $B_3(x;S/I)$ because $z=z-x_1^{p-1}x_4^p
  e_{123}=\partial(x_1^{p-1}x_4^{p-1}e_{1234})$, $x_1^{p-1}x_4^p$ being
  an element in $I$.}

\end{Remark}

Example \ref{bi} suggests the following:

\begin{Lemma}
\label{n-1} Let $I$ be an arbitrary monomial ideal. Then
$H_{n-1}(x;S/I)$ has a basis given by  cycles of length $\leq
\lfloor n/2 \rfloor$.
\end{Lemma}

\begin{proof}
Suppose $y$ is a multigraded $(n-1)$-cycle in the form given by
Lemma \ref{bound} (2). Then there exist $1\leq k_s<\ldots < k_1< n$
such that $y=\Sigma_{j=1}^s\gamma_ju_j e_{\sigma_j}$, $\gamma_j\in
K^*$, $u_j\not \in I$ monomials and
$\sigma_j=\{1,\ldots,n\}\setminus \{k_j\}$. Suppose that $y$ cannot
be written as a sum of cycles of  length $<s$. We will show that we
may choose a cycle $y'$ of length $\leq \lfloor n/2 \rfloor$ which
coincides with $y$ modulo $B_{n-1}(x;S/I)$. This is enough because
then these cycles will give a system of generators of
$H_{n-1}(x;S/I)$ from which we may choose a basis. We claim that
$\gamma_r=(-1)^{k_1-k_r} \gamma_1$ for all $1\leq r\leq s$. Let
$E=\{r: 1\leq r\leq s, \gamma_r=(-1)^{k_1-k_r}\gamma_1\}$. If $E\not
= \{1,\ldots, s\}$ then the element $q=\Sigma_{j\in E}\gamma_ju_j
e_{\sigma_j}$ of $K_{n-1}(x;S/I)$ is different from $y$. Thus $y$
cannot be a cycle because otherwise $y=q+(y-q)$ is a decomposition
of $y$ in a sum of two cycles of smaller length which is false
($1\in E$) by our assumption. So one of $\sigma_j$, $j\in E$ has a
neighbour  $\sigma_t$ in $z$ which is not in $q$, that is $t\not \in
E$, $x_{k_t}u_j\not \in I$ and
 $((-1)^{k_t-1}\gamma_j-(-1)^{k_j-1}\gamma_t)x_{k_t}u_j
e_{\sigma_j\setminus\{k_t\}}=0$ because $z$ is a cycle (actually the
above equation is written for the case $k_j>k_t$, otherwise all the
signs changed but the equation is not really affected). Thus
$\gamma_t=(-1)^{k_j-k_t}\gamma_j=(-1)^{k_1-k_t}\gamma_1$, that is
$t\in I$ which is false. Hence $E=\{1,\ldots,s\}$, that is our claim
holds.

As $y$ is multigraded note that $x_{k_j}|u_{k_j}$. We have the
following cycle
$$y'=\partial(\gamma_1( u_1/x_{k_1}) e_{1\ldots n})-y=\gamma_1\Sigma
(-1)^{k-1} (x_ku_1/x_{k_1}) e_{\sigma_k},$$ where
$\sigma_k=\{1,\ldots,n\}\setminus \{k\}$ and the sum is made over
all $k\in \{1,\ldots,n\}\setminus \{k_1,\ldots,k_s\}$. If $y'=0$
then $y\in B_{n-1}(x;S/I)$ and there exist nothing to show. If
$y'\not =0$ then length$(y)+$length$(y')\leq n$. Thus
$\min\{$length$(y)$, length$(y')\}\leq \lfloor n/2 \rfloor.$ As
$y\equiv -y'$ mod $B_{n-1}(x;S/I)$ we are done.

\end{proof}

\begin{Example}
\label{four} {\em Let $I=(x_3x_4x_5, x_2x_4x_5, x_1x_2x_4,
x_1x_2x_3, x_1x_3x_5)$ be an ideal in $S=K[x_1,\ldots,x_5]$ .
 We claim that there exists no monomial cycles or binomial cycles in the
 homology class of the following multigraded cycle
$$z=x_3x_4e_{125}-x_2x_4e_{135}+x_1x_3e_{245}.$$
We adopt the following notation: for a monomial element $ue_{abc}$
we will write let us say ${\bar a}$ , that is $ue_{{\bar a}bc}$  if
$x_au\not\in I$. So we may write
$$z=x_3x_4e_{{\bar 1}{\bar 2}5}-x_2x_4e_{1{\bar 3}5}+x_1x_3e_{2{\bar 4}5}$$
and now we can see easily that $z$ is indeed a cycle. We list all
monomial elements of $K_3(x;S/I)$, which are in the multigraded
class of $z$:
$$x_4x_5e_{{\bar 1}23}, x_3x_5e_{1{\bar 2}4}, x_3x_4e_{{\bar 1}{\bar 2}5},
x_2x_5e_{{\bar 1}{\bar 3}4}, x_2x_4e_{1{\bar 3}5}, x_2x_3e_{1{\bar
4}{\bar 5}}, x_1x_5e_{{\bar 2}3{\bar 4}},$$ $$ x_1x_4e_{2{\bar
3}{\bar 5}}, x_1x_3e_{2{\bar 4}5}, x_1x_2e_{34{\bar 5}}.$$ Clearly
no one is monomial cycle. A binomial cycle should have the form
$\gamma_1u_1e_{a{\bar b}c}+\gamma_2u_2e_{ac{\bar d}},$ $\gamma_i\in
K^*$, $u_i$ monomials. Thus we might find such pairs $(a{\bar
b}c),(ac{\bar d})$ among above. But there are no such pairs. For
example $e_{{\bar 1}23}$ could make such a pair only with $e_{{\bar
2}3{\bar 4}}$, $e_{2{\bar 3}{\bar 5}}$ but they are not of the
necessary type because each one has two numbers in bold. In this way
we see that our claim holds. It follows that $z$ cannot be written
modulo $B_3(x;S/I)$ as a sum of cycles of smaller length.

}
\end{Example}

\begin{Remark}
\label{char} {\em In \cite[Exercise 5.5.4]{BH} we see that for
monomial ideals $I$ the Betti numbers
$\beta_2,\beta_{n-2},\beta_{n-1}$ do not depend on the
characteristic of $K$. Since in the previous section we see that
$H_2(x;S/I)$ has monomial cycle basis we may ask by analogy if
$H_{n-2}(x;S/I)$ or $H_{n-1}(x;S/I)$ have monomial cycle basis.
Example \ref{bi}  shows that this is not true. On the other hand
note that in $K_{n-2}(x;S/I)$ there are cycles which cannot be
written modulo $B_{n-2}(x;S/I)$ as a sum of cycles of length $\leq
\lfloor n/2 \rfloor$ (see Example \ref{four}) as happens in the
$K_{n-1}(x;S/I)$ (see Lemma \ref{n-1}).}
\end{Remark}

We might ask which is the minimal possible positive integer $r$ such
that $H_3(x;S/I)$ has a basis given by 3-cycles of length $\leq r$.
The following example shows that $r$ could be even 4.

\begin{Example}
\label{five}{\em Let $I=(x_3x_4x_5x_6, x_2x_4x_5x_6, x_1x_2x_4x_6,
x_1x_2x_3x_4,x_1x_2x_3x_5, x_1x_3x_5x_6)$ be an ideal in
$S=K[x_1,\ldots,x_6]$ .
 We claim that there exists no monomial cycles or binomial cycles in the
 homology class of the following multigraded cycle
$$z=x_3x_4x_5e_{{\bar 1}{\bar 2}6}-x_2x_4x_5e_{{\bar 1}{\bar 3}6}+x_1x_3x_5e_{2{\bar 4}6}-x_1x_2x_4e_{3{\bar 5}6}$$
We list all monomial elements of $K_3(x;S/I)$, which are in the
multigraded class of $z$:
$$x_4x_5x_6e_{{\bar 1}23}, x_3x_5x_6e_{1{\bar 2}4}, x_3x_4x_6e_{{\bar 1}{\bar 2}5},
x_3x_4x_5e_{{\bar 1}{\bar 2}6},x_2x_5x_6e_{{\bar 1}{\bar 3}4},
x_2x_4x_6e_{1{\bar 3}5}, x_2x_4x_5e_{{\bar 1}{\bar 3}6},$$
$$x_2x_3x_6e_{{\bar 1}{\bar 4}{\bar 5}}, x_2x_3x_5e_{1{\bar
4}{\bar 6}},x_2x_3x_4e_{1{\bar 5}{\bar 6}},x_1x_5x_6e_{{\bar
2}3{\bar 4}}, x_1x_4x_6e_{2{\bar 3}{\bar 5}}, x_1x_4x_5e_{{\bar
2}{\bar 3}{\bar 6}},$$
 $$x_1x_3x_6e_{{\bar 2}{\bar 4}5},x_1x_3x_5e_{2{\bar 4}6},x_1x_3x_4e_{2{\bar
5}{\bar 6}}, x_1x_2x_6e_{{\bar 3}4{\bar 5}},x_1x_2x_5e_{3{\bar
4}{\bar 6}},x_1x_2x_4e_{3{\bar 5}6},x_1x_2x_3e_{45{\bar 6}}.$$
Clearly no one is monomial cycle. There are only 6 elements from the
above list (with just one "bar"), which can be used to construct
binomial cycles in the multigraded homology class of $z$. As in
Example \ref{four} we see that there are no binomial cycles. Now a
cycle of length 3 should have the form $\gamma_1u_1e_{a{\bar
b}c}+\gamma_2u_2e_{a{\bar c}{\bar d}}+\gamma_3u_3e_{ad{\bar t}},$
$\gamma_i\in K*$, $u_i$ monomials, for some $a,b,c,d,t\in
\{1,\ldots,n\}$, we may also have $t\in \{a,b\}$. We claim that this
is not possible. For example if $a=1$, $b=2$, $c=4$ then $d\in
\{3,5,6\}$. If $d=5$ note that in fact we have $e_{{\bar a}{\bar
c}{\bar d}}=e_{1{\bar 4}{\bar 5}}$ which is not possible to be term
in a cycle of length 3. If $d=6$ then $t\in \{2,3,5\}$. If for
example $t=5$ then we have $e_{a{\bar d}{\bar t}}=e_{1{\bar 5}{\bar
6}}$ which is not possible again. In this way we can show that our
claim holds. It follows that $z$ cannot be written modulo
$B_3(x;S/I)$ as a sum of cycles of smaller length.

}
\end{Example}

\section{Cycles  of Koszul homology modules of
principal $p$-Borel ideals}

Let $u=\Pi_{q=1}^n x_q^{\lambda_q}$ be a monomial and $J=\Pi_{q=1}^n
(x_1,\ldots,x_q)^{\alpha_q}$ for some integers
$\lambda_q,\alpha_q\geq 0$. Set $u_{\leq a}=\Pi_{q=1}^a
x_q^{\lambda_q}$, $J_{\leq a}=\Pi_{q=1}^a
(x_1,\ldots,x_q)^{\alpha_q}$ for some integer $1\leq a< n$ and $u_{>
a}=\Pi_{q>a}^n x_q^{\lambda_q}$, $J_{> a}=\Pi_{q>a}^n
(x_1,\ldots,x_q)^{\alpha_q}$.

\begin{Lemma}
\label{>}  Suppose that $ux_a\in J$ and $ux_{a+1}\not\in J$. Then
$u_{>a}\in J_{>a}$ and $u_{\leq a}\not \in J_{\leq a}$.
\end{Lemma}
\begin{proof}
By hypothesis we have $ux_a=vw$ for some monomials $v\in J_{\leq a}$
and $w\in J_{>a}$. If there exists $1\leq j\leq a$ with $x_j|w$ then
$$ux_{a+1}=(x_jv/x_a)(x_{a+1}w/x_j)\in J_{\leq a}J_{>a}=J.$$
Contradiction! It follows $w|u_{>a}=(ux_a)_{>a}$ and so $u_{>a}\in
J_{>a}$. If $u_{\leq a}\in J_{\leq a}$ then $u\in J$ and so
$ux_{a+1}\in J$ which is false.
\end{proof}

\begin{Lemma}
\label{crit} Let $v,w$ be some monomials in $(x_q)_{q>a}$ such that
$v|u$. Suppose that $wu/v\in J$ and $ux_a^r\in J$ for some integer
$r>0$. Then $ux_{a+1}^r\in J$.
\end{Lemma}
\begin{proof}
Suppose that $ux_{a+1}^r\not\in J$. We may suppose that there exists
an integer $0\leq d<r$ such that for $u'=ux_a^dx_{a+1}^{r-d-1}$ we
have $u'x_a\in J$ and $u'x_{a+1}\not \in J$.
 By Lemma \ref{>} we get $u'_{\leq
a}\not \in J_{\leq a}$. Since $wu/v=u_{\leq a}(wu_{>a}/v)\in
J\subset J_{\leq a}$ it follows $u_{\leq a}\in J_{\leq a}$ because
the variables from $wu_{>a}/v$ are regular on $S/J_{\leq a}$. Thus
$u'_{\leq a}\in J$ because $u'_{\leq a}=x_a^d u_{\leq a}$.
Contradiction!
\end{proof}

\begin{Lemma}
\label{crit2} Let $t>a$ be an integer and $v,w$ be some monomials in
$(x_q)_{q\geq t}$ such that $v|u$. Suppose that $wu/v\in J$ and
$ux_a^r\in J$ for some integer $r>0$. Then $ux_{t}^r\in J$.
\end{Lemma}

The proof goes applying Lemma \ref{crit} by recurrence.

Let $I=\Pi_{q=1}^n\Pi_{j=0}^s
((x_1,\ldots,x_q)^{[p^j]})^{\alpha_{qj}}$ be a principal $p$-Borel
ideal of $S$, $0\leq \alpha_{qj}<p$ being some integers. Let $G(I)$
be the set of minimal monomial generators of $I$. Note that all
monomials from $G(I)$ have the same degree.

\begin{Lemma}
\label{3-cycle1} Let $a,t,q,r$ be four positive integers such that
$a<t<r$, $a<q<r$, $q\not =t$ and $\gamma\in I$ a monomial which is a
multiple of $x_qx_r$.Suppose that $x_a\gamma/x_r\in I$ and
$x_t\gamma/x_q\in I$. Then either $x_t\gamma/x_r\in I$, or
$x_a\gamma/x_q\in I$.
\end{Lemma}
\begin{proof}
Apply induction on $s$. If $s=0$ then from $x_t\gamma/x_q\in I$ we
get $x_a\gamma/x_q\in I$, since $I$ is strongly stable in this case.

Now suppose $s>1$ and set $J=\Pi_{q=1}^n
((x_1,\ldots,x_q)^{[p^j]})^{\alpha_{q0}}$,

\noindent $T= \Pi_{q=1}^n\Pi_{j=1}^{s-1}
((x_1,\ldots,x_q)^{[p^j]})^{\alpha_{q,j+1}}$. Then $I=JT^{[p]}$. We
have $\gamma=\beta \alpha^p$ for some monomials $\beta\in J$,
$\alpha\in G(T)$.

Suppose from now on that $x_t\gamma/x_r\not \in I$. We will show
that $x_a\gamma/x_q\in I$. Then $x_r$ does not divide $\beta$ since
otherwise $x_t\gamma/x_r=(x_t\beta/x_r)\alpha^p\in I$, $J$ being
strongly stable. Contradiction! If $x_q|\beta$ then $x_a\beta/x_q\in
J$ because $J$ is strongly stable and so
$x_a\gamma/x_q=(x_a\beta/x_q)\alpha^p\in I$. Remains to study the
case when $x_q$ does not divide $\beta$. Then $x_q|\alpha$. We have
 $\alpha/x_q\not \in T$ because $\alpha\in G(T)$. From $x_t\gamma/x_q=x_t\beta
x_q^{p-1}(\alpha/x_q)^p\in I$ it follows either
\begin{enumerate}
\item[(i)] there exist $b, 1\leq b\leq n$ such that $x_b^p|\beta$,
$x_b\alpha/x_q\in T$ and $x_tx_q^{p-1}\beta/x_b^{p}\in T$, or
\item[(ii)] $x_t^{p-1}|\beta$, $x_q^{p-1}\beta/x_t^{p-1}\in J$ and
$x_t\alpha/x_q\in T$.
\end{enumerate}
Here we use the fact that all minimal monomial generators of $G(T)$
have the same degree and since $\alpha\in G(T)$ we get
$v\alpha/x_q\in T$  for $v$ being just one variable. If (i) holds
then $x_a\gamma/x_q=(x_ax_q^{p-1}\beta/x_b^p)(x_b\alpha/x_q)^p\in I$
since $J$ is strongly stable.

>From now on suppose that (ii) holds. As
$x_a\gamma/x_r=x_ax_r^{p-1}\beta (\alpha/x_r)^p\in I$ we get as
above the following two cases:
\begin{enumerate}
\item[(i')] $x_a^{p-1}|\beta$, $x_r^{p-1}\beta/x_a^{p-1}\in J$ and
$x_a\alpha/x_r\in T$, or
\item[(ii')] there exist $c, 1\leq c\leq n$ such that $x_c^p|\beta$,
$x_c\alpha/x_r\in T$ and $x_ax_r^{p-1}\beta/x_c^{p}\in T$.

\end{enumerate}

{\bf Case} (i') holds

It follows $x_r^{p-1}\beta/x_t^{p-1}\in J$ since $J$ is strongly
stable. Note that $\alpha$ satisfies over $T$ the condition of
$\gamma$ over $I$. By induction hypothesis we get either
$x_t\alpha/x_r\in T$ or $x_a\alpha/x_q\in T$. If $x_t\alpha/x_r\in
T$ then we get
$x_t\gamma/x_r=(x_r^{p-1}\beta/x_t^{p-1})(x_t\alpha/x_r)^p\in I$
which is false. Thus we must have $x_a\alpha/x_q\in T$. Then
$x_a\gamma/x_q=(x_q^{p-1}\beta/x_a^{p-1})(x_a\alpha/x_q)^p\in I$
since $x_r^{p-1}\beta/x_a^{p-1}\in J$, $J$ being strongly stable.

{\bf Case} (ii') holds

Suppose first $c\geq t$. From (ii),(ii') we obtain
$x_ax_r^{p-1}\beta/x_c^p\in J$ and $x_q^{p-1}\beta/x_t^{p-1}\in J$.
Using Lemma \ref{crit2} for $u=x_r^{p-1}\beta/x_c^p$,
$v=(x_rx_t)^{p-1}$, $w=(x_cx_q)^{p-1}$ we get
$ux_t=x_tx_r^{p-1}\beta/x_c^p\in J$. Then
$x_t\gamma/x_r=(x_tx_r^{p-1}\beta/x_c^p)(x_c\alpha/x_r)^p\in I$.
Contradiction!

Now suppose $c<t$. By (ii),(ii') we get $x_c\alpha/x_r\in T$,
$x_t\alpha/x_q\in T$. Apply the induction hypothesis for the
integers $c,t,q,r$ and the ideal $T$. It follows either
$x_c\alpha/x_q\in T$, or $x_t\alpha/x_r\in T$. In the first case we
get $x_a\gamma/x_q=(x_ax_q^{p-1}\beta/x_c^p)(x_c\alpha/x_q)^p\in I$
because $x_ax_r^{p-1}\beta/x_c^p\in J$, $J$ being strongly stable.

In the second case we can obtain
$x_t\gamma/x_r=(x_r^{p-1}\beta/x_t^{p-1})(x_t\alpha/x_r)^p\in I$,
that is a contradiction, providing we show that
$x_r^{p-1}\beta/x_t^{p-1}\in J$. But from
$x_ax_r^{p-1}\beta/x_c^p\in J$ we get
$x_ax_r^{p-1}\beta/x_cx_t^{p-1}\in J$ since $J$ is strongly stable.
Moreover if $c\leq a$ then we get even $x_r^{p-1}\beta/x_t^{p-1}\in
J$. If $c>a$ then we get the same thing using Lemma \ref{crit2} for
$u=x_r^{p-1}\beta/(x_cx_t^{p-1})$, $w=x_cx_t^{p-1}$, $v=x_r^{p-1}$
and $c$ instead $t$.
\end{proof}

\begin{Proposition}
\label{3} Let $z=\Sigma_{j=1}^s\gamma_ju_j e_{\sigma_j}$ be a
multigraded 3-cycle of $K_i(x;S/I)$, $\gamma_j\in K^*$, $u_j$
monomials, $\sigma_j\subset \{1,\ldots n\}$, $|\sigma_j|=i$ for all
$1\leq j\leq s$. Suppose that $s>1$, $m(u_j)\leq m(\sigma_j)$ for
all $j$, $\sigma_1=\{a,t,r\}$, $a<t<r$, $\sigma_1=\max_{1\leq j\leq
s}\sigma_j$ and $x_a u_1\in I$. Then there exists a multigraded
3-cycle $y $ of length $\leq 3$ such that $in (z)=in(y)$. Moreover
if the length of $y$ is 3 then the homology class of $y$ contains a
monomial cycle.
\end{Proposition}

\begin{proof}
By Lemmas \ref{bound} and \ref{xr} we have $m(\sigma_j)=r$ and $x_r
u_j\in I$ for all $1\leq j\leq s$. Set $\gamma=x_r u_1\in I$. Then
$x_a\gamma/x_r\in I$ by hypothesis. We may suppose $x_tu_1\not \in
I$ because otherwise $y=in(z)$ is a monomial cycle. Then $\sigma_1$
has a neighbour $\sigma_j$ in $z$ for $j>1$ by Lemma
\ref{neighbour}, let us say $\sigma_j=\{a,q,r\}$. As
$\sigma_j<\sigma_1$ we get $t<q$. We have $u_j=x_tu_1/x_q$ because
$z$ is multigraded. It follows $x_t\gamma/x_q=x_r u_j\in I$. By
Lemma \ref{3-cycle1} we get $x_a\gamma/x_q\in I$ since
$x_t\gamma/x_r=x_tu_1\not \in I$. Therefore $\partial((u_1/x_q)
e_{atqr})=y-(x_ru_1/x_q) e_{atq}$ for
$y=u_1e_{\sigma_1}-u_je_{\sigma_j}+(x_au_1/x_q) e_{tqr}$. From the
above we see that $(\gamma/x_q) e_{atq}$ is a monomial cycle and so
$y$ is a cycle.
\end{proof}

\begin{Remark}
\label{mightnot}{\em If $x_au_1\not \in I$ but $x_t u_1\in I$ then
$(x_r u_1/x_q) e_{atq}$ might be not a monomial cycle as happens in
the proof of the above proposition. Note that in the Example
\ref{bi} $(x_ru_1/x_q) e_{atq}=x_3^{p-1}x_4^p e_{123}$ is not a
monomial cycle because $x_3(x_3^{p-1}x_4^p)\not \in I$. However in
this example $z$ is already a binomial cycle, that is a cycle of
length $\leq 3$.}
\end{Remark}

\begin{Lemma}
\label{3-cycle2} Let $a,t,q,r$ be four positive integers such that
$a<t<r$, $a<q<r$ and $\gamma\in I$ a monomial which is a multiple of
$x_qx_r$. Suppose that $x_t\gamma/x_r\in I$, $x_a\gamma/x_q\in I$
and $x_t\gamma/x_q \not \in I$. Then the following statements hold:

\begin{enumerate}
\item[(i)] If $q>t$ then   $x_a\gamma/x_r\in I$.
\item[(ii)] If $q<t$ then $x_tx_a\gamma/x_rx_q\in I$.
\end{enumerate}
\end{Lemma}

\begin{proof}
Apply induction on $s$. If $s=0$ then from $x_t\gamma/x_r\in I$ we
get $x_a\gamma/x_r\in I$ and also $x_tx_a\gamma/x_rx_q\in I$  since
$I$ is strongly stable in this case.

Now suppose $s>1$ and set $J,T$ and $\gamma=\beta \alpha^p$ for some
monomials $\beta\in J$, $\alpha\in G(T)$ as in the proof of Lemma
\ref{3-cycle1}.
 Like there we may suppose that $x_r$ does not divide $\beta$.
 From $x_t\gamma/x_r=x_t\beta
x_r^{p-1}(\alpha/x_r)^p\in I$ it follows either
\begin{enumerate}
\item[(j)] there exist $b, 1\leq b\leq n$, $b\not =t$ such that $x_b^p|\beta$,
$x_b\alpha/x_r\in T$ and $x_tx_r^{p-1}\beta/x_b^{p}\in J$, or
\item[(jj)] $x_t^{p-1}|\beta$, $x_r^{p-1}\beta/x_t^{p-1}\in J$ and
$x_t\alpha/x_r\in T$.
\end{enumerate}

In case (j) we have $$ x_a\gamma/x_r=((x_a/x_t)(x_t\beta
x_r^{p-1}/x_b^p)) (x_b\alpha/x_r)^p\in I,$$ $J$ being strongly
stable.

 {\bf Case} (jj) when
$x_q|\beta$

  If $t<q$  then $x_t\beta/x_q\in J$ because $J$ is strongly
stable. Thus $x_t\gamma/x_q\in I$. Contradiction! If $t>q$ then we
get
$$x_tx_a\gamma/x_qx_r=((x_a/x_q)(x_r^{p-1}\beta/x_t^{p-1})(x_t\alpha/x_r)^p\in
I.$$

{\bf Case} (jj) when $x_q\not |\beta$.

 Then $x_q|\alpha$. As $\alpha/x_q\not \in T$ it follows from $x_a\gamma/x_q\in I$
either

\begin{enumerate}
\item[(j')]  there exist $c, 1\leq c\leq n$ such that $x_c^p|\beta$,
$x_c\alpha/x_q\in T$ and $x_ax_q^{p-1}\beta/x_c^{p}\in J$, or
\item[(jj')] $x_a^{p-1}|\beta$, $x_q^{p-1}\beta/x_a^{p-1}\in J$ and
$x_a\alpha/x_q\in T$.
\end{enumerate}

{\bf Subcase} (j'), $c\geq t$

Let $u=x_q^{p-1}\beta/x_c^p$. We have
$x_t\gamma/x_q=(x_tu)(x_c\alpha/x_q)^p\in I$ because by (jj) it
follows $x_tu=(x_a
x_t^{p-1}/x_c^p)(x_q^{p-1}/x_r^{p-1})(x_r^{p-1}\beta/x_t^{p-1})\in J
$,  $J$ being strongly stable. Contradiction!

 {\bf Subcase} (j'), $t>c$, $t<q$.

 Apply induction on $s$ for $\alpha, T$. Since $x_t\alpha/x_r\in
T$, $x_c\alpha/x_q \in T$ we get either $x_t\alpha/x_q\in T$ or
$x_c\alpha /x_r\in T$. If $x_t\alpha/x_q\in T$ then
$x_t\gamma/x_q=(x_q^{p-1}\beta/x_t^{p-1})(x_t\alpha/x_q)^p\in I$
because $x_r^{p-1}\beta/x_t^{p-1}\in J$, $J$ being strongly stable.
Contradiction! Note that  we did not use the condition $t<q$ in
order to get this contradiction from $x_t\alpha/x_q\in T$. Now
suppose that $x_c\alpha/x_r\in T$. Then we have
$x_a\gamma/x_r=(x_ax_r^{p-1}\beta/x_c^p)(x_c\alpha/x_r)^p\in I$ if
we show that $v=x_ax_r^{p-1}\beta/x_c^p\in J$. By (j') we have
$x_ax_q^{p-1}\beta/x_c^p\in J$. If $a\geq c$ then we get
$x_q^{p-1}\beta/x_c^{p-1}\in J$, $J$ being strongly stable. If $a<c$
we get the same thing applying Lemma \ref{crit2} to
$u'=x_q^{p-1}\beta/x_c^p$ because $x_au'\in J$ and
$u'(x_c^p/x_q^{p-1})=\beta\in J$.

Set ${\hat u}=x_r^{p-1}\beta/x_t^{p-1}x_c^{p-1}$. We have
$x_c^{p-1}{\hat u}\in J$ and ${\hat
u}(x_q^{p-1}x_t^{p-1}/x_r^{p-1})\in J$ from above. By Lemma
\ref{crit2} we get $x_t^{p-1}{\hat u}\in J$, that is
$x_r^{p-1}\beta/x_c^{p-1}\in J$. If $c\geq a$ it follows $v\in J$
since $J$ is strongly stable. If $c<a$ apply Lemma \ref{crit2} for
$u"=x_r^{p-1} \beta/x_c^p$ having $x_cu"\in J$ and
$u"(x_ax_q^{p-1}/x_r^{p-1})\in J$. It follows $v=x_au"\in J$.

{\bf Subcase} (j'), $t>c$, $t>q$.

 Apply induction on $s$ for $\alpha, T$. Since $x_t\alpha/x_r\in
T$, $x_c\alpha/x_q \in T$ we get either $x_t\alpha/x_q\in T$ or
$x_tx_c\alpha /(x_qx_r)\in T$. We saw above that $x_t\alpha/x_q\in
T$ gives a contradiction. Suppose that $x_tx_c\alpha/(x_qx_r)\in T$.
Then we have
$x_tx_a\gamma/(x_qx_r)=(x_ax_r^{p-1}x_q^{p-1}\beta/(x_t^{p-1}x_c^p))(x_tx_c\alpha/(x_rx_q))^p\in
I$ if we show that

\noindent $w=x_ax_r^{p-1}x_q^{p-1}\beta/(x_t^{p-1}x_c^p)\in J$. By
(j') we have $x_ax_q^{p-1}\beta/x_c^p\in J$. As in the previous case
we get $x_q^{p-1}\beta/x_c^{p-1}\in J$ if either $a\geq c$ or
$a<c\leq q$. If $a<q<c$ the same tricks bring only that $x_qu'\in
J$, that is $x_q^p\beta/x_c^p\in J$.

Suppose $a\geq c$ or $a<c\leq q$.
 We have $x_c^{p-1}{\hat u}\in J$ and ${\hat
u}(x_q^{p-1}x_t^{p-1}/x_r^{p-1})\in J$ as above. By Lemma
\ref{crit2} we get now $x_q^{p-1}{\hat u}\in J$, that is
$x_r^{p-1}x_q^{p-1}\beta/(x_t^{p-1}x_c^{p-1})\in J$. Therefore $w\in
J$ if $a\leq c\leq q$ because $J$ is strongly stable. If $a>c$ then
$y=x_r^{p-1}x_q^{p-1}\beta/(x_t^{p-1}x_c^p)$ satisfies $x_cy\in J$
and $yx_t^{p-1}x_a/x_r^{p-1}\in J$. Therefore by Lemma \ref{crit2}
we get $w\in J$. If $a<q<c$ then $z={\hat u}/x_c$ satisfies
$x_c^pz\in J$ and $z(x_q^px_t^{p-1}/x_r^{p-1})\in J$. By Lemma
\ref{crit2} we get as above
$x_q^px_r^{p-1}\beta/(x_t^{p-1}x_c^p)=x_q^pz\in J$ and so
 $w\in J$ because $J$ is strongly stable.

{\bf Subcase} (jj'), $t<q$.

As above $x_t\alpha/x_q\not\in T$. Apply induction on $s$ for
$\alpha, T$. Since $x_t\alpha/x_r\in T$, $x_c\alpha/x_q \in T$ we
get  $x_a\alpha /x_r\in T$.  It follows that
$$x_a\gamma/x_r=(x_r^{p-1}\beta/x_a^{p-1})(x_a\alpha/x_r)^p\in I$$
if $x_r^{p-1}\beta/x_a^{p-1}\in J$. Set ${\bar
u}=x_r^{p-1}\beta/x_a^{p-1}x_t^{p-1}$ and note that $x_a^{p-1}{\bar
u}\in J$ by (jj) and ${\bar u}(x_t^{p-1}x_q^{p-1}/x_r^{p-1})\in J$
by (jj'). By Lemma \ref{crit2} we get $x_t^{p-1}{\bar u}\in J$ which
is enough.

{\bf Subcase} (jj'), $t>q$.

As in the previous case we use induction hypothesis to get
$x_tx_a\alpha/(x_rx_q)\in T$. We have
$$x_tx_a\gamma/(x_rx_q)=(x_r^{p-1}x_q^{p-1}\beta/(x_t^{p-1}x_a^{p-1}))(x_tx_a\alpha/(x_rx_q))^p\in
I$$ if we show that
$v'=x_r^{p-1}x_q^{p-1}\beta/(x_t^{p-1}x_a^{p-1})\in J$. Apply Lemma
\ref{crit2} for

\noindent ${\tilde u}=x_r^{p-1}\beta/(x_t^{p-1}x_a^{p-1})$ because
$x_a^{p-1}{\tilde u}\in J$ and ${\tilde
u}(x_q^{p-1}x_t^{p-1}/x_r^{p-1})\in J$. We obtain $x_q^{p-1}{\tilde
u}\in J$ which is enough.

\end{proof}

\begin{Proposition}
\label{31} Let $z=\Sigma_{j=1}^s\gamma_ju_j e_{\sigma_j}$ be a
multigraded 3-cycle of $K_i(x;S/I)$, $\gamma_j\in K^*$, $u_j$
monomials, $\sigma_j\subset \{1,\ldots n\}$, $|\sigma_j|=i$ for all
$1\leq j\leq s$. Suppose that $s>1$, $m(u_j)\leq m(\sigma_j)$ for
all $j$, $\sigma_1=\{a,t,r\}$, $a<t<r$, $\sigma_1=\max_{1\leq j\leq
s}\sigma_j$ and $x_t u_1\in I$. Then there exists a multigraded
3-cycle $y $ of length $\leq 3$ such that $in (z)=in(y)$. Moreover
if the length of $y$ is 3 then the homology class of $y$ contains a
monomial cycle.
\end{Proposition}

\begin{proof}
We follow the proof of Proposition \ref{3}. Set $\gamma=x_r u_1\in
I$. Then $x_t\gamma/x_r\in I$ by hypothesis and we may suppose
$x_au_1\not \in I$. Thus $\sigma_1$ has a neighbour $\sigma_j$ in
$z$ for $j>1$ by Lemma \ref{neighbour}, let us say
$\sigma_j=\{t,q,r\}$. We have $u_j=x_au_1/x_q$ because $z$ is
multigraded. It follows $x_a\gamma/x_q=x_r u_j\in I$. Suppose $t<q$.
By Lemma \ref{3-cycle2} we get $x_a\gamma/x_r\in I$ since
$x_t\gamma/x_q=x_au_1\not \in I$. Therefore $\partial((u_1/x_q)
e_{atqr})=y-(x_ru_1/x_q) e_{atq}$ for
$y=u_1e_{\sigma_1}+u_je_{\sigma_j}-(x_tu_1/x_q) e_{tqr}$. From the
above we see that $(\gamma/x_q) e_{atq}$ is a monomial cycle and so
$y$ is a cycle. Now suppose $q<t$. Then the same lemma gives that
$x_tu_j\in I$, that is $u_1e_{\sigma_1}-u_je_{\sigma_j}$ is a
binomial cycle.
\end{proof}

\begin{Proposition}
\label{32} Let $z=\Sigma_{j=1}^s\gamma_ju_j e_{\sigma_j}$ be a
multigraded 3-cycle of $K_i(x;S/I)$, $\gamma_j\in K^*$, $u_j$
monomials, $\sigma_j\subset \{1,\ldots n\}$, $|\sigma_j|=i$ for all
$1\leq j\leq s$. Suppose that $s>1$, $m(u_j)\leq m(\sigma_j)$ for
all $j$, $\sigma_1=\{a,t,r\}$, $a<t<r$ and $\sigma_1=\max_{1\leq
j\leq s}\sigma_j$. Then there exists a multigraded 3-cycle $y $ of
length $\leq 4$ such that $in (z)=in(y)$. Moreover if the length of
$y$ is 3 then the homology class of $y$ contains a monomial cycle
and if the length is 4 then  the homology class of $y$ contains a
binomial cycle.
\end{Proposition}

\begin{proof}
We follow the proof of Propositions \ref{3} and  \ref{31}. Set
$\gamma=x_r u_1\in I$. Using the quoted propositions we may suppose
$x_tu_1\not\in I$ and $x_au_1\not \in I$. Thus $\sigma_1$ has two
neighbours $\sigma_j$, $\sigma_k$ in $z$ for $j,k>1$ by Lemma
\ref{neighbour}, let us say $\sigma_j=\{a,q,r\}$ and
$\sigma_k=\{c,t,r\}$.  We have $u_k=x_au_1/x_c$, $u_j=x_tu_1/x_q$
 because $z$ is multigraded. It follows $x_a\gamma/x_c=x_r u_k\in I$ and
 $x_t\gamma/x_q\in I$ by Lemma \ref{xr}. As $u_1e_{\sigma_1}=in (z)$
 we get $\sigma_j<\sigma_1$ and so
$t<q$.

{\bf Case} $c=q$.

Thus we have   $x_ax_ru_1/x_q=x_a\gamma/x_q\in I$. We see that
$y'=(x_ru_1/x_q)e_{atq}$ is a monomial cycle and so $y=y'+
\partial((u_1/x_q)e_{atqr})$ is a cycle of length $\leq 3$ such that
$in (z)=in (y)$.

 {\bf Case} $a<t<\min\{c,q\}$

Suppose $c<q$. By Lemma \ref{3-cycle2} we get either
$x_a\gamma/x_q\in I$ or $x_t\gamma/x_c \in I$. If $c>q$ then by
Lemma \ref{3-cycle1} we get the same thing. Above we already studied
the case when $x_ax_ru_1/x_q\in I$.  If $x_tx_ru_1/x_c=x_t\gamma/x_c
\in I$ then similarly $y"=(x_ru_1/x_c)e_{atc}$ is monomial cycle and
$y=y"+
\partial((u_1/x_q)e_{atcr})$ is a cycle of length $\leq 3$ such that
$in( z)=in( y)$.

{\bf Case} $a<c<t<q$.

Using Lemma \ref{3-cycle2} as above we get either $x_t\gamma/x_c \in
I$ (case already studied above) or $x_tx_a\gamma/(x_qx_c)\in I$. In
the second case it follows that
$\phi=(x_ru_1/x_q)e_{atq}-(x_ax_ru_1/(x_qx_c))e_{ctq}$ is a binomial
cycle. Then  $$y_1=\phi+ \partial((u_1/x_q)
e_{atqr})=u_1e_{\sigma_1}-u_je_{\sigma_j}+(x_au_1/x_q)e_{tqr}-(x_ax_ru_1/(x_qx_c))e_{ctq}$$
is a cycle and the cycle
$$y=y_1-\partial((x_au_1/(x_cx_q))e_{ctqr}=u_1e_{\sigma_1}-u_je_{\sigma_j}-u_ke_{\sigma_k}+
(x_ax_tu_1/(x_cx_q))e_{cqr}$$ is of length 4. But as in
 Example \ref{four} we may see that
 $$y-\partial( (u_1/x_q)e_{atqr})+\partial( (x_au_1/(x_cx_q))e_{ctqr}=-(x_rx_au_1/(x_cx_q))e_{ctq}
 +(x_ru_1/x_q)e_{atq}.$$

\end{proof}

\begin{Theorem}
\label{main1}  $H_3(x;S/I)$ has a basis of binomial cycles.
\end{Theorem}

For the proof apply Proposition \ref{32}.

\section{Monomial cycle basis on Koszul homology modules of some
principal $p$-Borel ideals.}

The principal $p$-Borel ideals $I\subset S$ such that $S/I$ is
Cohen-Macaulay have the form $I=\Pi_{j=0}^s(m^{[p^j]})^{\alpha_j}$,
$0\leq \alpha_j<p$. For these ideals  is well known the description
of a canonical monomial cycle basis of $H_i(x;S/I)$. Fix $2\leq
i\leq n$. Let $0\leq t\leq s$ be an integer and for $v\in
G(m^{\alpha_t})$ denote $v'=v/x_{m(v)}$. Let $B_{it}(I)$ be the
following set of elements from $K_i(x;S/I)$
$$\{wv'^{p^t} x_{\sigma}^{p^t-1} e_{\sigma} : \ w\in  G(\Pi_{j>t}(m{[p^j]})^{\alpha_j}), v\in G(m^{\alpha_t}),
\sigma\subset {\bar n}, |\sigma|=i, m(\sigma)=m(v) \}$$ and
$B_i(I)=\cup_{t=0}^s B_{it}(I)$.

\begin{Theorem}[Aramova-Herzog \cite{AH}]
\label{ah} The elements of $B_i(I)$ are cycles in $K_i(x;S/I)$ and
their homology classes form a basis in $H_i(x;S/I)$ for $i\geq 2$.
\end{Theorem}

\begin{Remark}
\label{obstr}{\em This result holds independently of the
characteristic of $K$ as we had pointed the definition of $p$-Borel
ideals is pure combinatorial. But note that Theorem \ref{ah} does
not hold if $\alpha_j\geq p$ for some $j$. Indeed, the ideal
$I=(x_1,x_2)^4\subset S=K[x_1,x_2]$ is strongly stable and a
monomial basis of $H_2(x;S/I)$ is given by $T=\{x_1^3 e_1\wedge e_2,
x_1^2x_2 e_1\wedge e_2, x_1x_2^2 e_1\wedge e_2, x_2^3 e_1\wedge
e_2\}$ by \cite{AH1} (see also  \cite{EK}). Since
$I=(x_1,x_2)^2(x_1^2,x_2^2)$ one can compute $B_0(I)=T$ and
$B_1(I)=\{x_1x_2 e_1\wedge e_2\}$ but $x_1x_2 e_1\wedge e_2$ is not
cycle in $K_2(x;S/I)$. So the condition $\alpha_j<p$ is necessary
and this is an obstruction for an extension of Theorem \ref{ah}. }
\end{Remark}

The question appeared  in Remark \ref{obstr} perhaps can be solved
extending somehow the Theorem \ref{ah} for the case when $\alpha_j$
are arbitrary. In some special cases a possible tool could be the
following lemma.

\begin{Lemma}
\label{cas} Let  $I=\Pi_{j=0}^s(m^{[p^j]})^{\alpha_j}$, where
$\alpha_j\geq 0$ are arbitrary integers. If $n=2$ then there exist
some integers $0\leq j_0<j_1<\ldots <j_k$ and some positive integers
$(\gamma_t)_{0\leq t\leq k}$ such that $\gamma_t<p^{j_{t+1}-j_t}$
for $t<k$ and $I=\Pi_{t=0}^k(m^{[p^{j_t}]})^{\gamma_t}.$

\end{Lemma}

For the proof apply by recurrence the relation
$m^{p^t}m^{[p^t]}=(m^{p^t})^2$.

Set $I=\Pi_{t=0}^k(m^{[p^{j_t}]})^{\gamma_t}$ as above but for any
$n$ and let $C_{it}(I)$ be the following set of elements from
$K_i(x;S/I)$
$$\{wv'^{p^{j_t}} x_{\sigma}^{p^{j_t}-1} e_{\sigma} : \ w\in  G(\Pi_{r>t}(m^{[p^{j_r}]})^{\gamma_r}), v\in G(m^{\gamma_{t}}),
\sigma\subset {\bar n}, |\sigma|=i, m(\sigma)=m(v) \}$$ and
$C_i(I)=\cup_{t=0}^s C_{it}(I)$. A  variant of Theorem \ref{ah} is
the following theorem:

\begin{Theorem}
\label{ah1} The elements of $C_i(I)$ are cycles in $K_i(x;S/I)$ and
their homology classes form a basis in $H_i(x;S/I)$ for $i\geq 2$.
\end{Theorem}

Since Lemma \ref{cas} works only in the case $n=2$ this  gives
almost nothing more than \ref{ah}. Unfortunately the ideals of type
$T=m^{p^j}m^{[p^j]}$ could be  bad for example when $p=3$ and $n=3$
then $T=m^6\setminus \{x_1^2x_2^2x_3^2\}$.

 Let $M$ be a graded $S$-module and
$\beta_{ij}(M)=\dim_K \Tor^i_S(K,M)_j$ the $ij$-th graded Betti
number of $M$.

\begin{Corollary}
\label{indep}$\beta_{ij}(S/I)$ does not depend on the characteristic
of the field $K$ for all $i,j$.
\end{Corollary}

For the proof note that $H_i(x;S/I)\iso \Tor_i^S(K,S/I)$ and so
$\beta_{ij}(S/I)$ is the sum of some $|C_{it}(I)|$ which has nothing
to do with the characteristic of $K$.

\begin{Remark}
\label{obstrBetti}{\em Note that $\beta_{ij}(S/I)$ does not depend
on the characteristic of $K$ when $I$ is stable by \cite{EK}. In
\cite{Po} it shows that the extremal graded Betti numbers of $S/I$
(see \cite{BCP}) do not depend on the characteristic of $K$ when $I$
is a Borel type ideal (see \cite{HPV}). In particular this happens
for $p$-Borel ideals and so we might ask if all $\beta_{ij}(S/I)$ do
not depend on the characteristic of $K$ in the case of  $p$-Borel
ideals. The Corollary \ref{indep} is a small hope.}
\end{Remark}

>From now on let $I$ be the $p$-Borel ideal generated by the monomial
$x_{n-1}^{\gamma}x_n^{\alpha}$ for some integer $\gamma,\alpha\geq
0$, that is
$I=\Pi_{j=0}^s((m_{n-1}^{[p^j]})^{\gamma_j}(m^{[p^j]})^{\alpha_j})$,
where $m_{n-1}=(x_1,\ldots,x_{n-1})$, and $\gamma_j,\alpha_j$ are
defined by the $p$-adic expansion of $\gamma$, respectively
$\alpha$. The main result of this section is the following:

\begin{Theorem}
\label{main} Suppose that $\alpha_j+\gamma_j<p$ for all $0\leq j\leq
s$. Then \begin{enumerate}
\item{} $H_i(x;S/I)$ has a monomial cycle basis for all $i\geq 2$,
and
\item{} $\beta_{ij}(S/I)$ does not depend on the characteristic of
$K$ for all $i,j$.
\end{enumerate}
\end{Theorem}

For the proof we need some preparations. Suppose $\alpha>0$. Let
$r=\max\{j: \alpha_j>0\}$ and set
$J=\Pi_{j=0}^s(m_{n-1}^{[p^j]})^{\gamma_j}$ and
$$I'=(\Pi_{j=0}^{r-1}(m^{[p^j]})^{\alpha_j})(m^{[p^r]})^{\alpha_r-1}.$$
Then

\begin{Lemma}
\label{div} $(I:x_n^{p^r})=JI'$.
\end{Lemma}

\begin{proof} Obviously if $L,T$ are some monomial ideals and $v$ is
a monomial then $(LT:v)=(L:v)T+L(T:v)$. Applying this fact we get
$$(I:x_n^{p^r})=J\Sigma\
\Pi_{j=0}^r\Pi_{k=1}^{\alpha_j}(m^{[p^j]}:x_n^{c_{jk}})=J\Sigma\
\Pi_{j=0}^r\Pi_{k=1}^{\alpha_j}(m_{n-1}^{[p^j]},x_n^{p^j-c_{jk}}),$$
where the sum is taken over all integers $0\leq c_{jk}\leq p^j$ such
that $\Sigma_{j=0}^r\Sigma_{k=1}^{\alpha_j} c_{jk}=p^r$. For each
$j$ let $\Lambda_j\subset \{k: 1\leq k\leq \alpha_j, c_{jk}>0\}$ be
any subset. Set $u=\Sigma_{j=0}^r\Sigma_{k\in \Lambda_j}
(p^j-c_{jk})$. We claim that
$x_n^u\Pi_{j=0}^r(m_{n-1}^{[p^j]})^{\alpha_j-|\Lambda_j |}\subset
I'$. Clearly if our claim holds then $(I:x_n^{p^r})\subset JI'$, the
other inclusion being trivial. Note that the claim holds because
$u\geq (\Sigma_{j=0}^r  |\Lambda_j| p^j)-p^r$.

\end{proof}

Let $a$ be an integer such that $0\leq a\leq\alpha$ and
$a=\Sigma_{j=0}^r a_jp^j$, $0\leq a_j<p$ the $p$-adic expansion of
$a$. Set $\alpha_a=\Sigma_{j, \alpha_j\geq a_j} (\alpha_j-a_j)p^j$
and $\alpha_{aj}=\alpha_j-a_j$ if $\alpha_j\geq a_j$ and 0
otherwise. Set
$$I_a=J(\Pi_{j=0}^r(m^{[p^j]})^{\alpha_{aj}}),$$
where $J$ is defined above. Let $\pi:S\to {\bar
S}=K[x_1,\ldots,x_{n-1}]$ be the ${\bar S}$-morphism given by
$x_n\to 0$.

\begin{Lemma}
\label{pdiv} $\pi(I:x_nâ)$ is the $p$-Borel ideal generated by the
monomial $x_{n-1}^{\gamma} x_{n-1}^{\alpha_a}$, that is
$\pi(I:x_n^a)=\pi(I_a)$.
\end{Lemma}

\begin{proof}
It is enough to show the above equality for the case $\gamma=0$. As
in the proof of Lemma \ref{div} we have
$$(I:x_n^a)=\Sigma\
\Pi_{j=0}^r\Pi_{k=1}^{\alpha_j}(m_{n-1}^{[p^j]},x_n^{p^j-c_{jk}}),$$
where the sum is taken over all integers $0\leq c_{jk}\leq p^j$ such
that $\Sigma_{j=0}^r\Sigma_{k=1}^{\alpha_j} c_{jk}=a$. Note that
$\pi(m_{n-1}^{[p^j]},x_n^{p^j-c_{jk}})$ is $m_{n-1}^{[p^j]}$ if
$c_{jk}<p^j$ and ${\bar S}$ otherwise. It follows that
$\pi(m_{n-1}^{[p^j]},x_n^{p^j-c_{jk}})\subset \pi(I_a)$, the
equality holds only when  $\min\{a_j,\alpha_j\}=|\{k:c_{jk}=p^j\}|$
for all $j,k$. Hence $\pi(I:x_n^a)=\pi(I_a)$.
\end{proof}

Let $T\subset S$ be an arbitrary ideal and ${\bar T}=\pi(T)$.

\begin{Lemma}
\label{h} $H_i(x;{\bar S}/{\bar T})\iso H_i(x_1,\ldots,x_{n-1};
{\bar S}/{\bar T})\oplus H_{i-1}(x_1,\ldots,x_{n-1}; {\bar S}/{\bar
T})$ and in particular $\beta_{ij}^S({\bar S}/{\bar
T})=\beta_{ij}^{\bar S}({\bar S}/{\bar T})+\beta_{i-1,j}^{\bar
S}({\bar S}/{\bar T})$, where $\beta_{ij}^S({\bar S}/{\bar T})$ is
the $i,j$-th graded Betti number of ${\bar S}/{\bar T}$ over $S$.
\end{Lemma}

\begin{proof}
By \cite[Proposition 1.6.21]{BH} we have $$H_i(x;{\bar S}/{\bar
T})\iso H_i(x; S/({\bar T},x_n))\iso H_i(x_1,\ldots,x_{n-1};
S/({\bar T}S))\otimes_S(\wedge^1S)\iso$$ $$
H_{i}(x_1,\ldots,x_{n-1}; {\bar S}/{\bar T})\otimes_S(\wedge^1S),$$
the last isomorphism follows because $S$ is flat over ${\bar S}$.
This is enough because $\beta_{ij}^{\bar S}({\bar S}/{\bar
T})=\dim_k\Tor_i^{\bar S}(K,{\bar S}/{\bar T})_j=\dim_k
H_i(x_1,\ldots,x_{n-1}; {\bar S}/{\bar T})_j$.

\end{proof}

Because of the above isomorphism we may write $$H_i(x;{\bar S}/{\bar
T})= H_i(x_1,\ldots,x_{n-1}; {\bar S}/{\bar T})\oplus
H_{i-1}(x_1,\ldots,x_{n-1}; {\bar S}/{\bar T})\wedge e_n.$$ where by
abus of notation we write $$H_{i-1}(x_1,\ldots,x_{n-1}; {\bar
S}/{\bar T})\wedge e_n$$ for $\{cls(z\wedge e_n): z\ cycle\  of\
K_{i-1}(x_1,\ldots,x_{n-1}; {\bar S}/{\bar T})\}$.  We have the
following multigraded exact sequence $$(*)\ \ \ \ \ 0\to
S/(T:x_n)(-1)\to S/T\to {\bar S}/{\bar T}\to 0,$$ where first map is
given by multiplication with $x_n$. Applying Koszul homology long
exact sequence to (*) we get the following multigraded exact
sequence:
$$(**) H_i(x;S/(T:x_n)(-1))\to H_i(x;S/T)\to H_i(x;{\bar S}/{\bar
T})\to H_{i-1}(x;S/(T:x_n)(-1)),$$ where we denote by $\delta_i$ the
last map. Next lemma describes how acts $\delta_i$.

\begin{Lemma}
\label{del}  $\delta_i$ maps $H_i(x_1,\ldots,x_{n-1}; {\bar S}/{\bar
T})$ in zero and if $z$ is a cycle of

\noindent $K_{i-1}(x_1,\ldots,x_{n-1}; {\bar S}/{\bar T})$ then
$\delta_i$ maps $cls(z\wedge e_n)$ in
$$(-1)^{i-1}cls(z)\in
H_{i-1}(x_1,\ldots,x_{n-1}; S/( T:x_n)(-1)).$$

\end{Lemma}
\begin{proof} We have the following commutative diagram

\[
\begin{array}{ccccccccc}
0 & \to & K_i(x;S/(T:x_n)(-1))& \to & K_i(x;S/T) & \to &K_i(x;{\bar
S}/{\bar T}) &\to & 0\\
& & \downarrow & & \downarrow & &\downarrow &&\\
0 & \to & K_{i-1}(x;S/(T:x_n)(-1))& \to & K_{i-1}(x;S/T) & \to
&K_{i-1}(x;{\bar
S}/{\bar T}) &\to & 0\\
\end{array}
\]

Let $w$ be a cycle of $K_i(x;{\bar S}/{\bar T})$. By construction of
$\delta_i$ we must lift $w$ to an element $v\in K_i(x;S/T)$. Then
$\partial(v)=x_n y$ for a cycle $y\in K_{i-1}(x;S/(T:x_n)(-1))$ and
we may write $\delta_i(cls(w))=cls(y)$. Here we may take $v=w \in
K_i(x;S/T)$ which is a cycle. Then we have $y=0$ and so
$\delta_i(w)=0$. Now we take $w=z\wedge e_n$. As in the first case
we may take $v=z\wedge e_n$ but this time this is not cycle in
$K_i(x;S/T)$. We have $\partial(z\wedge e_n)=\partial(z)\wedge e_n+
(-1)^{i-1}x_nz=(-1)^{i-1}x_nz$ since $\partial(z)=0$. Then
$\delta_i(cls(z\wedge e_n))=(-1)^{i-1}cls(z)$.

\end{proof}

Let $f_i$ be the composite map $ K_i(x;S/T) \to K_i(x;{\bar S}/{\bar
T})\to  H_i(x_1,\ldots,x_{n-1}; {\bar S}/{\bar T})$ where the second
map $q_2$ is the second projection of the direct sum given by Lemma
\ref{h}. Then $f_i$ has a canonical section $\rho_i^T$ given by
$cls(z)\to cls(z)\in  H_i(x;S/T)$, $z$ being a cycle of
 $K_i(x_1,\ldots,x_{n-1}; {\bar S}/{\bar T})$. Let $\eta_i^T:H_i(x_1,\ldots,x_{n-1}; {\bar S}/{\bar T})
 \to H_i(x;{\bar S}/\pi(T:x_n)$
be the canonical map associated to the surjection ${\bar S}/{\bar T}
 \to {\bar S}/\pi(T:x_n)$.

\begin{Corollary}
\label{delta} The following statements hold:
\begin{enumerate}
\item{} $\delta_{i+1}=(-1)^i \rho_i^{(T:x_n)}\eta_i^T q_2$,
\item{} $\Ker \delta_{i+1}\iso H_i(x_1,\ldots,x_{n-1}; {\bar S}/{\bar
T})\oplus \Ker \eta_i^T$,
\item{} $\Im \delta_{i+1}\iso \Im \eta_i^T$.
\end{enumerate}
\end{Corollary}

\begin{Lemma}
\label{lift} Let $ue_{\sigma}$, $u\in S$ monomial with $m(u)<n$.
Suppose that $ue_{\sigma}$ is a monomial cycle of
$K_i(x_1,\ldots,x_{n-1};{\bar S}/{\bar T})$ and induces an element
of $\Ker \eta_i^T$, that is $ue_{\sigma}=\partial(z)$ in
$K_i(x_1,\ldots,x_{n-1};{\bar S}/\pi( T:x_n) )$ for an element

\noindent $z\in K_{i+1}(x_1,\ldots,x_{n-1};{\bar S}/\pi( T:x_n) )$.
Then $ue_{\sigma}\wedge e_n+(-1)^ix_nz$ is a cycle in
$K_{i+1}(x;S/T)$. If $ue_{\sigma}$ is zero in
$K_i(x_1,\ldots,x_{n-1};{\bar S}/\pi( T:x_n) )$ then
$ue_{\sigma}\wedge e_n$ is a cycle also in $K_{i+1}(x;S/T)$.
\end{Lemma}

\begin{proof}
We have $\partial(ue_{\sigma}\wedge
e_n+(-1)^ix_nz)=\partial(ue_{\sigma})\wedge e_n+
(-1)^{i-1}x_nz+(-1)^ix_n\partial(z)= \partial(ue_{\sigma})\wedge
e_n=0$ because $ue_{\sigma}$ is a monomial cycle in
$K_i(x_1,\ldots,x_{n-1};{\bar S}/{\bar T})$. The second statement
holds because then $z=0$.
\end{proof}

 Now we may return to
give the proof of Theorem \ref{main}.

\begin{proof} Apply induction on $c=\Sigma_{j=0}^s \alpha_jp^j$. If
$c=0$ then we are in the case of Theorem \ref{ah} and Corollary
\ref{indep}. Suppose $c>0$ and set $r=\max \{j:\alpha_j\not =0\}$.
By Lemma \ref{div} $I'=(I:x_n^{p^{r}})$ is the $p$-Borel ideal
generated by the monomial $x_{n-1}^{\gamma}x_n^{\alpha-p^{r}}$ and
from induction hypothesis $H_i(x;S/I')$ has monomial cyclic basis
and $\beta_{ij}^S(S/I')$ does not depend on the characteristic of
$K$ for all $i,j$.

Let $0\leq a\leq p^{r}$ be an integer. By decreasing induction we
show that $H_i(x;S/(I:x_n^a)$ has monomial cycle basis  and
$\beta_{ij}^S(S/(I:x_n^a))$ does not depend on the characteristic of
$K$ for all $i,j$. Above we saw the case $a=p^{r}$. Suppose
$a<p^{r}$. The exact multigraded sequence (**) given before Lemma
\ref{del} with Corollary \ref{delta} give for $T=(I:x_n^a)$ the
following exact multigraded sequence
$$0\to \Im\delta_{i+1}\iso\Im \eta_i^{(I:x_n^a)}\to
H_i(x;S/(I:x_n^{a+1}))\to H_i(x;S/(I:x_n^a))\to \Ker \delta_i\iso$$
$$H_i(x_1,\ldots, x_{n-1};{\bar S}/\pi(I:x_n^a))\oplus \Ker
\eta_{i-1}^{(I:x_n^a)}\to 0,$$ where
$\eta_i^{(I:x_n^a)}:H_i(x_1,\ldots, x_{n-1};{\bar
S}/\pi(I:x_n^a))\to H_i(x_1,\ldots, x_{n-1};{\bar
S}/\pi(I:x_n^{a+1}))$ is given in Corollary \ref{delta}. By Lemma
\ref{pdiv} we see that $\pi(I:x_n^a)=\pi(I_a)$ is the $p$-Borel
ideal generated by a power of $x_{n-1}$ and it is subject to Theorem
\ref{ah} and Corollary \ref{indep} because $\gamma_j+\alpha_{aj}\leq
\gamma_j+\alpha_j<p$ for all $j$. In particular,
$B_i(\pi(I:x_n^a))=B_i(\pi(I_a))$ is in $H_i(x_1,\ldots,
x_{n-1};{\bar S}/\pi(I:x_n^a))$ a monomial cycle basis.

Using the induction hypothesis on $a$ we see that
$H_i(x;S/(I:x_n^{a+1}))$ has a monomial cyclic basis and
$\beta_{ij}^S(S/(I:x_n^{a+1}))$ does not depend on the
characteristic of $K$ for all $i,j$. Then the conclusion follows
from the above multigraded exact sequence and Lemma \ref{lift} if we
show the following statements:
\begin{enumerate}
\item{} A monomial cycle basis of $H_i(x_1,\ldots, x_{n-1};{\bar
S}/\pi(I:x_n^a))$ can be lifted to a monomial cycle subset of
$H_i(x;S/(I:x_n^a))$.
\item{} $\Ker \eta_i^{(I:x_n^a)}$ has a monomial cycle basis which
can be lifted to a monomial cycle subset of $H_i(x;S/(I:x_n^a))$.
\item{} $\dim_K\Ker \eta_i^{(I:x_n^a)}$ does not depend on the
characteristic of $K$.

\end{enumerate}

Actually  (1) was already seen in the proof of Lemma \ref{del}. For
(2) and (3) we study how act $\eta_i^{\pi(I_a)}$ on $B_i(\pi(I_a))$.
We have the following cases:

{\bf Case} $a_0<p-1$.

Then $a+1=(\Sigma_{j>0} a_jp^j)+(a_0+1)$ is the $p$-adic expansion
of $a+1$, that is $(a+1)_j=a_j$ for all $j>0$ and $(a+1)_0=a_0+1$.
We have $\alpha_{a+1,j}=\alpha_{aj}$ for $j>0$. If $\alpha_0\leq
a_0$ then $\alpha_0=\alpha_{a+1,0}=0$ and $\eta_i^{\pi(I_a)}$ acts
identical because $\pi(I_a)=\pi(I_{a+1})$. Thus
$\Ker\eta_i^{\pi(I_a)}=0$. If $\alpha_0>a_0$ then
$\alpha_{a+1,0}=\alpha_{a,0}-1$. Thus $\eta_i^{\pi(I_a)}$ acts
identical on $\cup_{t\geq 1}B_{it}(\pi(I_a))$ and send
$B_{io}(\pi(I_a))$ in zero since if $v\in G(m_{n-1}^{\alpha_{a0}})$
then $v'\in G(m_{n-1}^{\alpha_{a0}-1})$. So the monomial cyclic
basis of $\Ker\eta_i^{\pi(I_a)}$ is given by $B_{io}(\pi(I_a))$.

{\bf Case} $a_j=p-1$ for $0\leq j<t$, $a_t<p-1$.

Then $a+1=(a_t+1)p^t+\Sigma_{j>t}a_jp^j$ is the $p$-adic expansion
of $a+1$, that is $(a+1)_j=0$ for $j<t$, $(a+1)_j=a_t+1$ for $j=t$
and $(a+1)_j=a_j$ for $j>t$. We have $\alpha_{a+1,j}=\alpha_{a,j}$
for $j>t$ and so $\eta_i^{\pi(I_a)}$  acts identical on $\cup_{j>
t}B_{ij}(\pi(I_a))$. If $\alpha_t\leq a_t$ then
$\alpha_{at}=\alpha_{a+1,t}$ and $\eta_i^{\pi(I_a)}$ acts identical
on $B_{it}(\pi(I_a))$. If $\alpha_t>a_t$ then
$\alpha_{a+1,t}=\alpha_{at}-1$ and $\eta_i^{\pi(I_a)}$  send
$B_{it}(\pi(I_a))$ to zero. Suppose $j<t$. Then
$\alpha_{a,j}=\alpha_{a+1,j}$ and $\eta_i^{\pi(I_a)}$ acts identical
on $B_{ij}(\pi(I_a))$. Otherwise we have
$\alpha_{a,j}>\alpha_{a+1,j}$ and $\eta_i^{\pi(I_a)}$ send
$B_{ij}(\pi(I_a))$ to zero.

Consequently, given $j\geq 0$ in both cases $\eta_i^{\pi(I_a)}$
either acts identical on $B_{ij}(\pi(I_a))$ or send it to zero. It
follows that $\Ker\eta_i^{\pi(I_a)}$ has a monomial cyclic basis
which can be lifted to $H_i(x;S/(I:x_n^a))$ by Lemma \ref{lift}. It
consists of some $B_{ij}(\pi(I_a))$  whose cardinal does not depend
on the characteristic of $K$. This ends our decreasing induction.
Thus the ideal $(I:x_n^a) $ satisfy the conditions (1), (2) from the
Theorem \ref{main} for all $0\leq a\leq p^{r}$. In particular this
holds for $a=0$.

\end{proof}

\begin{Remark}
\label{nonpr} {\em Note that the above proof shows also that some
non-principal $p$-Borel ideals of the form $(I:x_n^a) $ have
monomial cyclic bases.}
\end{Remark}

We end this section with an example illustrating the proof of
Theorem \ref{main}.

\begin{Example}
\label{ill}{\em Let $n=3$, $p=2$, $S=K[x_1,x_2,x_3]$,
$m=(x_1,x_2,x_3)$, $I=m^{[2]}m$. Using Theorem \ref{ah} a cyclic
basis of $H_2(x;S/I)$ is given by $B_{21}(I)=\{x_1x_2e_1\wedge
e_2,x_1x_3e_1\wedge e_3, x_2x_3e_2\wedge e_3\}$ and
$B_{20}(I)=\{x_i^2e_1\wedge e_2,x_i^2e_1\wedge e_3, x_i^2e_2\wedge
e_3: 1\leq i\leq 3\}$. We will show this independently using the
procedure from the proof of Theorem \ref{main}. Let $\pi:S\to
{\bar S}=K[x_1,x_2]$ be the ${\bar S}$-morphism given by $x_3\to
0$. Then ${\bar I}=\pi(I)=m_2^{[2]}m_2$, where $m_2=(x_1,x_2)$ and
$\pi(I:x_3)=m_2^{[2]}$, $(I:x_3^2)=m$. Note that  monomial cyclic
basis of $H_2(x_1,x_2;{\bar S}/m_2^{[2]}m_2)$, $H_2(x_1,x_2;{\bar
S}/m_2^{[2]})$ are given by $B_{21}({\bar I})=\{x_1x_2e_1\wedge
e_2\}$ and $B_{20}({\bar I})=\{x_1^2e_1\wedge e_2, x_2^2e_1\wedge
e_2\}$ respectively $\{x_1x_2e_1\wedge e_2\}$. The map $\eta_2^I$
maps $B_{20}({\bar I})$ in zero and it is identity on
$B_{21}({\bar I})$. The maps $\eta_1^I$, $\eta_i^{(I:x_3)}$,
$i=1,2$ are zero maps.

We have the following multigraded exact sequence
$$\Im \eta_2^{(I:x_3)}=0\to H_2(x;S/m)\to H_2(x;S/(I:x_3))\to$$
$$H_2(x_1,x_2;{\bar S}/m_2^{[2]}) \oplus H_1(x_1,x_2;{\bar
S}/m_2^{[2]})\wedge e_3\to 0.$$ As the monomial cyclic basis of
$H_2(x;S/m)$, $H_2(x_1,x_2;{\bar S}/m_2^{[2]})$, $H_1(x_1,x_2;{\bar
S}/m_2^{[2]})$ are $\{e_1\wedge e_2, e_1\wedge e_3, e_2\wedge
e_3\}$, respectively $\{x_1x_2e_1\wedge e_2\}$, respectively
$\{x_1e_1,x_2e_2\}$ we see that a monomial cycle basis in
$H_2(x;S/(I:x_3))$ is given by $$T_2(I:x_3)=\{x_3 e_1\wedge e_2,
x_3e_1\wedge e_3, x_3 e_2\wedge e_3, x_1x_2e_1\wedge e_2,
x_1e_1\wedge e_3, x_2 e_2\wedge e_3\}.$$

Now consider the multigraded exact sequence
$$0\to \Im \eta_2^{I}\to H_2(x;S/(I:x_3))\to H_2(x;S/I)\to$$
$$H_2(x_1,x_2;{\bar S}/m_2^{[2]}m_2) \oplus H_1(x_1,x_2;{\bar
S}/m_2^{[2]}m_2)\wedge e_3\to 0.$$ As the monomial cyclic basis of
$\Im \eta_2^I$, $H_2(x_1,x_2;{\bar S}/m_2^{[2]}m_2)$,
$H_1(x_1,x_2;{\bar S}/m_2^{[2]}m_2)$ are $\{x_1x_2e_1\wedge e_2\}$,
respectively $T_2(m_2^{[2]})=\{x_1x_2e_1\wedge e_2,x_1^2e_1\wedge
e_2, x_2^2e_1\wedge e_2\}$, respectively
$T_1(m_2^{[2]})=\{x_1^2e_1\wedge e_3,x_2^2e_1\wedge e_3,
x_1^2e_2\wedge e_3,x_2^2e_2\wedge e_3\}$ we see that a monomial
cycle basis in $H_2(x;S/I)$ is given by
$$x_3[T_2(I:x_3)\setminus \{x_1x_2 e_1\wedge e_2\}]\cup T_2(m_2^{[2]})\cup
T_1(m_2^{[2]})=B_{20}(I)\cup B_{21}(I)=B_2(I).$$}
\end{Example}


\begin{thebibliography}{1}

\bibitem{AH1} A.\ Aramova, \ J.\ Herzog, Koszul cycles and
Eliahou-Kervaire type resolutions, J.Alg, {\bf 181} (1996),347-370.

\bibitem{AH} A.\ Aramova,\ J.\ Herzog, $p$-Borel principal ideals,
Illinois J. Math., {\bf 41-1}, (1997), 103-121.


 \bibitem{BCP} D.\ Bayer,\ H.\ Charalambous, S.\ Popescu, Extremal
 Betti
numbers and applications to monomial  ideals, J. Alg., {\bf
221}(2)(1999), 497-512.

\bibitem{BM} D.\ Bayer, D.\ Mumford, What can be computed in algebraic
geometry? in {\sl Comutational Algebraic Geometry and Commutative
Algebra }, Symposia Mathematica XXXIV (1993), 1-48.





 \bibitem{BS} D.\ Bayer, M.\ Stillman, A criterion for detecting $m$
 regularity, Invent. Math. {\bf 87} (1987), 1-11.


\bibitem{BH}  W.\ Bruns,\ J. Herzog, {\sl Cohen-Macaulay rings}, Revised
Edition,
 Cambridge, 1996.

\bibitem{CS} G.\ Caviglia, E.\ Sbarra, Characteristic-free bounds for
the Castelnuovo-Mumford regularity, to appear in Compositio Math.


\bibitem{Ei} D.\  Eisenbud , {\sl Commutative algebra , with a view toward
geometry}, Graduate Texts Math. Springer, 1995.



\bibitem{EK} S.\ Eliahou, M.\ Kervaire, Minimal resolutions of some
monomial ideals, J. Alg., {\bf 129} (1990), 1-25.

\bibitem{HP} J.\ Herzog, D.\ Popescu, On the regularity of $p$-Borel
ideals, Proceed. of AMS, {\bf 129-9}, 2563-2570.


\bibitem{HPV} J.\ Herzog, D.\ Popescu, M.\ Vladoiu, On the Ext-modules of
ideals of Borel type, Contemporary Math. {\bf 331} (2003),
171-186.

\bibitem{MM} E.\ Mayr, A.\ Meyer, The complexity of the word problem
for commutative semigroups and polynomial ideals, Adv. in Math.,
{\bf 46} (1982). 305-329.


\bibitem{Pa} K.\ Pardue, {\it Nonstandard Borel fixed ideals}, Dissertation,
Brandeis University, 1994.

\bibitem{Po} D.\ Popescu, Extremal Betti numbers and regularity of
Borel type ideals, Bull. Math. Soc. Sc. Math. Roumanie, {\bf
48(96)}, 1(2005),65-72.

\end{thebibliography}
\end{document}